\documentclass[12pt]{article}
\usepackage{amssymb}
\usepackage{latexsym}
\usepackage{fullpage}
\usepackage{amscd}
\def\R{\mathbb{R}}
\def\C{\mathbb{C}}
\def\H{\mathbb{H}}
\def\N{\mathbb{N}}

\def\sp{\mathfrak{sp}}
\def\SO{\mathbf{SO}}
\def\SU{\mathbf{SU}}
\def\Sp{\mathbf{Sp}}
\def\U{\mathbf{U}}
\def\Cas{\mathrm{Cas}}
\def\End{\mathrm{End}}
\def\Gr{\mathbf{Gr}}
\def\Hom{\mathrm{Hom}\,}
\def\Iso{\mathrm{Iso}}
\def\id{\mathrm{id}\,}

\def\ind{\mathrm{index}\,}
\def\Ric{\mathrm{Ric}}
\def\Td{\mathrm{td}\,}

\def\pr{\mathrm{pr}\,}
\def\#{\sharp}

\def\A{\mathcal{A}}
\def\E{\mathcal{E}}

\def\G{\mathbf{G}}
\def\g{\mathfrak{g}}

\def\j{\mathrm{jet}}
\def\k{\mathfrak{k}}
\def\L{\Lambda}
\def\l{\lambda}
\def\p{\mathfrak{p}}
\def\S{\mathrm{Sym}}
\def\t{\mathfrak{t}}
\def\T{\textstyle\bigotimes}
\def\W{\mathrm{Wolf}\,}

\def\<#1,#2>{\langle\,#1,\,#2\,\rangle}
\newtheorem{Lemma}{Lemma}[section]
\newtheorem{Proposition}[Lemma]{Proposition}
\newtheorem{Theorem}[Lemma]{Theorem}
\newtheorem{Corollary}[Lemma]{Corollary}
\newtheorem{Definition}[Lemma]{Definition}
\def\proof{\noindent\textbf{Proof:}\quad}
\def\summary#1{(\textsl{#1})\hfill\break}

\def\pfill{\par\vskip3mm plus1mm minus1mm\noindent}
\def\qed{\ensuremath{\quad\Box\quad}}
\begin{document}
\title{An upper bound for a Hilbert polynomial on quaternionic K{\"a}hler 
       manifolds}
\author{Uwe Semmelmann \& Gregor Weingart
 \thanks{The first author is a member of the {\sl European Differential
 Geometry Endeavour} (EDGE), Research Training Network HPRN-CT-2000-00101,
 supported by The European Human Potential Programme. The second author is
 partially supported by the SFB 611 \textsl{``Singul\"are Ph\"anomene
  und Skalierung in mathematischen Modellen''}}}
\maketitle
\begin{abstract}
 In this article we prove an upper bound for a Hilbert polynomial
 on quaternionic K{\"a}hler manifolds of positive scalar curvature.
 As corollaries we obtain bounds on the quaternionic volume and 
 the degree of the associated twistor space. Moreover the article
 contains some details on differential equations of finite type.
 Part of it is used in the proof of the main theorem.
\end{abstract}
\vskip0.5cm
%
%
\section{Introduction}
 A quaternionic K{\"a}hler manifold is a $4n$--dimensional Riemannian
 manifold $(M,\,g)$  with holonomy contained in $\Sp(1)\cdot\Sp(n)$. 
 Associated with $M$ is the twistor space $Z$ which is the total space
 of a $\C P^1$-fibration over $M$. If the scalar curvature of $M$ is
 positive, the twistor space $Z$ is a K{\"a}hler--Einstein manifold
 admitting a complex contact structure with values in a holomorphic
 line bundle $L$, in particular $L^{n+1}\,\cong\,K^*$ is isomorphic
 to the anticanonical bundle $K^*\,:=\,\L^{2n+1}T^{1,0}Z$ of $Z$.
 In this situation S.~Salamon introduced the Hilbert polynomial $P(r)$
 of $M$ in \cite{sala1}, which for integers $r$ is defined as the holomorphic
 Euler characteristic of $L^r$. It is easy to show that $P(r)$ for
 $r\in \N$ is given by the index of a certain twisted Dirac operator
 on $M$. The Hilbert polynomial $P(r)$ contains interesting information
 on the quaternionic K{\"a}hler manifold $M$, e.~g.~the dimension of
 the isometry group of $M$ is given as $P(1)$ according to a result
 of S.~Salamon \cite{sala1}.

 In \cite{betti} we developed a method based on representation theory of
 $\Sp(1)\cdot\Sp(n)$ for determining the kernel of twisted Dirac operators
 in terms of minimal eigenspaces of certain natural second order differential
 operators. This method is used here to prove our main theorem:

 \begin{Theorem}\label{main}
  Let $(M^{4n},\,g)$ be a quaternionic K{\"a}hler manifold of positive
  scalar curvature. The Hilbert polynomial $P(r)$ of $M$ is bounded from
  above by the Hilbert polynomial $P_{\H P^n}$ of the quaternionic
  projective space $\H P^n$ in the sense that for all integers $r\geq 0$:
  $$
   0 \;\; \leq \;\; P(r) \;\;\leq\;\; P_{\H P^n}(r)
   \;\;=\;\; { 2n+1 + 2r \choose 2n+1 } \ .
  $$
 \end{Theorem}

 As applications we obtain in Section~\ref{estimate} an upper bound for
 the dimension of isometry group, the degree of the twistor space and the
 quaternionic volume. Section~\ref{hwolf} contains a formula for the Hilbert
 polynomial on the Wolf spaces together with explicit expressions in several
 examples. The final Section~\ref{killing} is more or less independent from
 the rest of the article. It contains details on differential operators of
 finite type, which gives the background for the proof of our main theorem
 in Section~\ref{finito}.
\section{The Hilbert polynomial}\label{polynomial}
 Let $(M^{4n},\,g)$ be a quaternionic K{\"a}hler manifold, a Riemannian
 manifold with holonomy contained in $\Sp(1)\cdot \Sp(n) \subset\SO(4n)$.
 The holonomy reduction associates a globally defined vector bundle on $M$
 to every representation of $\Sp(1)\cdot \Sp(n)$. These representations are
 (sums of) tensor products of $\Sp(1)$-- and $\Sp(n)$--representations which
 factor through the projection $\Sp(1) \times \Sp(n) \rightarrow \Sp(1)\cdot
 \Sp(n)$. In particular the standard representations $H$ and $E$ of $\Sp(1)$
 and $\Sp(n)$ respectively induce only locally defined vector bundles whereas
 $\S^2H$ or $H \otimes E \cong TM\otimes_\R\C$ are globally defined.

 The twistor space $Z$ of $M$ is defined as the unit sphere bundle in
 $\S^2H \subset \End\,TM$. If $M$ is a quaternionic K{\"a}hler manifold
 of positive scalar curvature its twistor space $Z$ is a K{\"a}hler--Einstein
 manifold of complex dimension $2n+1$ endowed with a complex contact structure
 $\eta:\;T^{1,0}Z\,\longrightarrow\,L$ with values in a holomorphic line
 bundle $L$. Hence $L$ is an $n+1$--th root of the anticanonical bundle
 $K^*\,:=\,\L^{2n+1}T^{1,0}Z$:
 $$
  L^{n+1}\;\cong \; K^* \ .
 $$
 The {\it Hilbert polynomial} $P(r)$ of the polarized variety $(Z, L)$ is
 defined as the holomorphic Euler characteristic of the line bundle $L^r$
 $$
   P(r) \;:=\;  \chi(Z,\,{\mathcal O}(L^r) )
   \;= \; \sum^{2n+1}_{s=0}\,(-1)^s h^s(L^r)
   \;=\; \<e^{rc_1(L)}\,\Td Z,\,[Z]>\,,
 $$
 which can be calculated from the Todd class $\Td Z$ of $Z$ according to
 the Riemann--Roch Theorem. From the latter formula we conclude that $P(r)$
 is a polynomial of degree $\leq 2n+1$ in $r$ whose leading coefficient is
 given by the {\it degree} $\,\deg\,Z\,:=\,\<c_1(L)^{2n+1}, [Z]>$ of the
 twistor space $Z$:
 \begin{equation}\label{polynom1}
  P(r) \;=\; \frac{\deg(Z)}{(2n+1)!} \,r^{2n+1}\;+\;
  \textrm{terms of lower order}\,.
 \end{equation}
 Its constant term equals $P(0)\,=\,\<\Td Z,[Z]>\,=\,1$, because $Z$ is a
 compact K\"ahler manifold of positive scalar curvature. In addition $c_1(L)$
 is a positive class in $H^{1,1}(Z;\,\R)$, hence $L^r$ is a negative line
 bundle for $r < 0 $ and the Kodaira Vanishing Theorem asserts $h^s(L^r)
 \,=\,0$ for $r<0$ and $s \le 2n$. Moreover Serre duality implies
 $$
  h^s(L^r) \;=\;h^{2n+1-s}(L^{-r} \otimes K) \;=\; h^{2n+1-s}(L^{-r-n-1}),
 $$
 thus it follows the symmetry $P(r)\,=\,-P(-r-n-1)$.
 These properties were proved by S.~Salamon in \cite{sala1}, where he also showed 
 that on any quaternionic K{\"a}hler
 manifold of non vanishing scalar curvature the space $H^0(\,Z,\,
 \mathcal{O}(L)\,)$ of holomorphic sections of $L$ is isomorphic to the space
 of infinitesimal isometries of $M$ of dimension $P(1)\,=\,\dim \Iso(M,g)$.
 More generally an integration along the fibres results in the formula
 \begin{equation} \label{index1}
  P(r) \;=\; \<{\hat A}(M)\,\mathrm{ch}\,\S^{n+2r} H,\,[M]>\,,
 \end{equation}
 which identifies $P(r)$ with the index of a twisted Dirac operator on $M$.
 This observation will be our starting point for proving the upper bound for
 the Hilbert polynomial in Section \ref{estimate}. In order to provide a closer
 link between the polynomial and the geometry of $M$ we need to calculate
 the Chern character of $\S^{n+2r}H$:

 \begin{Lemma}\label{bernoulli}
  The Chern character $\,\mathrm{ch}\,\S^{n+2r}H$ of $\,\S^{n+2r}H$
  can be written
  $$
   \mathrm{ch}\,\S^{n+2r}H\;=\;
   \sum_{l\ge 0}\, \frac{2^{2l+1}}{(2l+1)!}\,B_{2l+1}(r+\frac{n+2}{2})\,u^l
  $$
  where $u\,:=\,p_1(H)$ and $B_{2l+1}(x)$ is the $2l+1$--th
  Bernoulli polynomial.
 \end{Lemma}

 \proof
 By the splitting principle we may think of the quaternionic bundle
 $H$ as a sum $\ell\oplus\ell^{-1}$ of two conjugated line bundles
 with first Chern classes $\pm\sqrt{u}$ or $u\,:=\,p_1(H)\,=\,c_1(\ell)^2$.
 Similarly we may think of $\S^kH,\,k\geq 0,$ as decomposed into $\ell^k
 \oplus \ell^{k-2}\oplus\ldots\oplus\ell^{-k}$. Using the characteristic
 property $B_{\mu+1}(w+1)\,-\,B_{\mu+1}(w)\,=\,(\mu+1)\,w^\mu$ of the
 Bernoulli polynomials together with $B_{\mu+1}(w)\,=\,(-1)^{\mu+1}\,
 B_{\mu+1}(1-w)$ the Chern character of $\S^kH$ becomes:
 \begin{eqnarray*}
  \mathrm{ch}\,\S^kH
  &=& e^{k\sqrt u}\,+\,e^{(k-2)\sqrt u}\,+\,\cdots\,+\,e^{-k\sqrt u}\\
  &=& \sum_{\mu\geq 0}\,\frac{2^\mu}{\mu!}\,\left(\,\sum_{\nu=0}^k\,
      (\frac k2-\nu)^\mu\,\right)\,\sqrt{u\,}^\mu\\
  &=& \sum_{\mu\geq 0}\,\frac{2^\mu}{(\mu+1)!}\,(\,B_{\mu+1}(\frac k2+1)
      \,-\,B_{\mu+1}(-\frac k2)\,)\sqrt{u\,}^\mu\\
  &=& \sum_{l\geq 0}\,\frac{2^{2l+1}}{(2l+1)!}
      \,B_{2l+1}(\frac k2+1)\,u^l\,.\qed
 \end{eqnarray*}

 Using the expansion $\,B_k(x)\,=\,x^k\,-\,\frac12{k\choose1}x^{k-1}\,+\,
 \frac16{k\choose2}x^{k-2}\,-\,\frac1{30}{k\choose4}x^{k-4}\,\pm
 \ldots$ of the Bernoulli polynomials we can expand formula
 (\ref{index1}) in powers of $r$ and interpret its coefficients
 in terms of geometric data of $M$ e.~g.~the leading coefficient
 is proportional to the volume of $M$. In fact the integral class
 $$
  4u\;\;=\;\;p_1(\S^2H)\;\;=\;\;(\frac{\kappa}{8\pi n(n+2)})^2\,\Omega
 $$
 is represented by a multiple of the Kraines form $\Omega$, whose
 $n$--th power $\Omega^n\,=\,(2n+1)!\,\mathrm{vol}_M$ is essentially
 the Riemannian volume form of $M$ for the canonical quaternionic
 orientation. Defining the {\it quaternionic volume} of $M$ as the integer
 $$
  v(M)\;\;:=\;\;\<(4u)^n,[M]>\;\;=\;\;(2n+1)!\,
  (\frac{\kappa}{8\pi n(n+2)})^{2n}\,\mathrm{vol}(M)
 $$
 and substituting the expansion of $\,\mathrm{ch}\,\S^{n+2r}H\,$ into
 equation (\ref{index1}) we find:
 \begin{eqnarray*}
 P(r) 
  &=&
  \sum_{l=0}^n\, \frac{2}{(2l+1)!}\,B_{2l+1}( r+\frac{n}{2}+1 )
  \,\<{\hat A}(M)(4u)^l, [M] > \\[1.5ex]
  &=&
  \frac{2}{(2n+1)!}\,B_{2n+1}(r+\frac{n}{2}+1)\,v(M)
  \;+\; \ldots \\[1.5ex]
  &=&
  \frac{2}{(2n+1)!}\, v(M)\,r^{2n+1}
  \;+\; \frac{n+1}{(2n)!}\,v(M)\, r^{2n}\;+\;\ldots \ .
 \end{eqnarray*}
 In particular we obtain the well--known formula $\,\deg(Z)=2v(M)$
 (c.~f.~\cite{sala2}). 
 Since $P(r)$ is integer valued for all integers $r$ it follows that $P(r)$ 
 can be written as a linear combination of binomial coefficients,~i.~e.~there 
 exist integers $n_i$ for $i=0,\ldots,2n+1$ with
 $$
  P(r) \;=\; \sum^{2n+1}_{i=0}\,n_i\,{r \choose i} .
 $$
 On the quaternionic projective space $\H P^n$ the bundle $L$ is the square
 of a globally defined holomorphic line bundle $L^{\frac12}$. Consequently
 the polynomial  $P(r)$ has additional zeroes for  $r=-1/2, \ldots, -n/2$
 which lead to the explicit formula:
 $$
  P(r) \;\;=\;\; { 2n+1 + 2r \choose 2n+1 } \qquad \mbox{for} \quad r \ge 0 \ .
 $$
\section{Indices of Twisted Dirac Operators}\label{estimate}
 According to a result of S.~Salamon mentioned in Section~\ref{polynomial} the
 values $P(r)$ of the Hilbert polynomial of a quaternionic K\"ahler manifold
 $M$ are indices of twisted Dirac operators $D_{\S^{n+2r}H}$. These twisted
 Dirac operators belong to a two parameter family of twisted Dirac operators
 of particular interest in quaternionic geometry and quite a lot is known
 about their indices. In \cite{betti} we proved a general principle for
 operators in this family identifying their kernels with minimal eigenspaces
 for certain selfadjoint second order differential operators $\Delta_\pi$.
 Applying this principle for the twisted Dirac operators $D_{\S^{n+2r}H}$
 allows us to prove the estimate on the Hilbert polynomial through a detailed
 study of these minimal eigenspaces in the next section.

 In \cite{betti} we defined for any representation $\pi$ of $\Sp(1)\cdot
 \Sp(n)$ a natural second order differential operator $\Delta_\pi$ acting
 on sections of the associated bundle $\pi(M)$ by
 \begin{equation}\label{Deltapi}
  \Delta_\pi \;:=\; \nabla^*\nabla \;+\; 2q(R) \ ,
 \end{equation}
 where $q(R)$ is a selfadjoint endomorphism of $\pi(M)$ depending linearly
 on $R$. In terms of a local orthonormal base $\{\omega_i\}$ of $\sp(1)\oplus
 \sp(n)\,\subset\,\L^2TM$ this curvature endomorphism can be written
 $2q(R)\,=\,\sum \omega_i\cdot R(\omega_i)$ where $\omega_i$ acts via
 the differential of the representation of $\Sp(1)\cdot\Sp(n)$ on $\pi$. 
 For a parallel subbundle $\pi(M)\,\subset\,\Lambda^\bullet T^*M$ the
 operator $\Delta_\pi$ coincides with the Hodge--Laplacian on  forms,
 since definition (\ref{Deltapi}) is nothing else but the classical
 Weitzenb\"ock formula. In particular we have $2q(R)\,=\,\Ric$ on 1--forms.
 Moreover it is not difficult to show that on symmetric spaces $G/K$ the
 operator $\Delta_\pi$ is the Casimir operator of $G$. For the twisted Dirac
 operators in question the results of \cite{betti} readily imply that the
 index of $D_{\S^{n+2r}H}$ is given by the dimension of the minimal eigenspace
 of the operator $\Delta_{\S^{2r}H}$ acting on sections of $\S^{2r}H$:

 \begin{Lemma}\label{lowerbound}
  Let $(M^{4n},\,g)$ be a quaternionic K{\"a}hler manifold of positive
  scalar curvature $\kappa$ and let $\Delta_{\S^{2r}H}$ be the differential
  operator defined above acting on sections of $\S^{2r}H$. The spectrum of
  $\Delta_{\S^{2r}H}$ is bounded below by $\l_{2r}\,:=\,\frac{\kappa}
  {2n(n+2)}r(n+1+r)$ and:
  $$
   \ker(D_{\S^{n+2r}H}) \;\;\cong\;\; \ker(\Delta_{\S^{2r}H} - \l_{2r}) \ .
  $$
 \end{Lemma}

 \proof
 First we recall some notations and general formulas of \cite{betti}.
 The spinor bundle of $M$ is associated to the representation
 $S\,=\,\oplus^n_{i=0} R^{l,n-l}$ with $R^{l,d}\,:=\,\S^lH \otimes \L^d_0E$.
 Quaternionic K{\"a}hler manifolds $M$ not isometric to $\H P^n$ are spin if
 and only if the quaternionic dimension $n$ is even, nevertheless this does
 not cause any problems as we will only consider twisted Dirac operators
 acting on globally defined vector bundles. If $\pi \subset S\otimes R^{l,d}$
 is any representation occurring in the decomposition of the tensor product
 into irreducible summands, then:
 $$
  \Delta_\pi \;=\; D^2_{R^{l,d}}\Big\vert_{\pi(M)} \;+\; \phi(l,d)
  \quad \mbox{with} \quad
  \phi(l,d):= {\kappa\over 8n(n+2)}(l+d-n)(l-d+n+2)\,.
 $$
 We call the representation  $R^{l,d}$ a maximal twist for $\pi$ if the
 number $\phi(l,d)$ is maximal among all representations $R^{\tilde l,\,
 \tilde d}$ with $\pi \subset S \otimes R^{\tilde l,\,\tilde d}$. Using
 this notion we have an identification \cite{betti}
 \begin{equation}\label{twist}
  \ker(D_{R^{l,d}}) \;\;\cong\;\; \bigoplus_\pi \ker(\Delta_{\pi} - \phi(l,d)) 
 \end{equation}
 where the sum is over all $\pi\subset S \otimes R^{l,d}$ for which $R^{l,d}$
 is a maximal twist. Indeed if $R^{l,d}$ is not a maximal twist for $\pi$ and
 so $\pi \subset S \otimes  R^{\tilde l,\,\tilde d}$ with $\phi(\tilde l,\,
 \tilde d) > \phi(l,\,d)$, then
 $$
  D^2_{R^{l,d}}\Big\vert_{\pi(M)} 
   \;=\; D^2_{R^{\tilde l,\tilde d}}\Big\vert_{\pi(M)}\;+\;
  (\;\phi(\tilde l,\,\tilde d)\,-\,\phi(l,\,d)\;)\,.
 $$
 Thus $D^2_{R^{l,d}}$ restricted to sections of $\pi(M)$ is positive and
 the representation $\pi$ cannot contribute to the kernel of $D_{R^{l,d}}$.
 Specializing to our case $R^{l,d}\,=\,\S^{n+2r}H$ we observe that the only
 representation $\pi$ with maximal twist $R^{n+2r,0}$ is $\pi\,=\,\S^{2r}H$,
 hence the lemma follows from formula (\ref{twist}).
 \qed

 \pfill
 For the canonical quaternionic orientation of $H\otimes E$ induced by the
 Kraines form $\Omega$, the half spin representations are given by:
  $$
   S^+\quad:=\quad\bigoplus_{r\equiv n\,(2)}\; R^{r,n-r}
   \qquad\qquad
   S^-\quad:=\quad\bigoplus_{r\not\equiv n\,(2)}\;R^{r,n-r}\,.
  $$
 Using the Glebsch--Gordan formula it follows immediately that
 $\S^{2r}H \subset S^+\otimes\S^{n+2r}H$. Hence the index of the
 twisted Dirac operator $D_{\S^{n+2r}H}$ is just the dimension of
 its kernel:
 \begin{equation}\label{index}
   P(r) \;=\; \ind(D_{\S^{n+2r}H}) \;=\;
  \dim \ker (\Delta_{\S^{2r}H} - \l_{2r}) \ .
 \end{equation}

 \pfill
 The following proposition contains an estimate which is our most important 
 technical result, which immediately implies Theorem~\ref{main}: 

 \begin{Proposition}\label{estimateprop}
  Let $M^{4n}$ be a quaternionic K{\"a}hler manifold of positive
  scalar curvature 
  $$
   \dim \ker (\Delta_{\S^{2r}H} - \l_{2r}) \;\le\; \dim
   \S^{2r}(\,H\oplus E\,)\;=\; { 2n+1 + 2r \choose 2n+1 } 
  $$
 \end{Proposition}

 The proof of this proposition will be given in Section \ref{finito}.
 As a first application we obtain the well--known upper bound for the
 dimension of the isometry group of a quaternionic K{\"a}hler manifold:

 \begin{Corollary}
  Let $(M^{4n},\,g)$ be a quaternionic K{\"a}hler manifold of positive
  scalar curvature:
  $$
   \dim \Iso(M,\,g) \;\leq\; \dim \Sp(n+1) \;=\; (n+1)(2n+3)
  $$
 \end{Corollary}

 In Section \ref{polynomial} we defined the degree of the twistor
 space $Z$ of a quaternionic K{\"a}hler manifold. Since it appears in
 the leading coefficient of the Hilbert polynomial and since the
 twistor space of $\H P^n$ is $\C P^{2n+1}$ our estimate immediately
 implies $\deg(Z) \le \deg(\C P^{2n+1})$. By the definition of degree
 this can be reformulated into an inequality of the corresponding
 Chern numbers:

 \begin{Corollary}
  Let $Z$ be the twistor space of a quaternionic K{\"a}hler
  manifold of positive scalar curvature:
  $$
   c_1(Z)^{2n+1} \;\leq\;
   c_1(\C P^{2n+1})^{2n+1}
   \;=\; 2^{2n+1}\,(n+1)^{2n+1}\ .
  $$
 \end{Corollary}

 For a compact K{\"a}hler--Einstein manifold $M$ of complex dimension
 $m$ and positive scalar curvature C.~LeBrun and S.~Salamon proved the
 following estimate for the top power of the first Chern class
 (c.f.~\cite{lebrun1})
 \begin{equation}\label{degree}
  c_1(M)^m \;\le\; \frac{m+1}{q}\,c_1(\C P^m)^m
 \end{equation}
 where $q$ is the largest integer dividing $c_1(M)$, the so called
 index of $M$. In case $M$ admits a complex contact structure $m\,=\,
 2n+1$ is odd and the index is given by $q=n+1$, so the estimate
 (\ref{degree}) becomes $c_1(M)^m \leq 2 c_1(\C P^m)^m$. However the compact
 K\"ahler--Einstein manifolds with complex contact structures are
 precisely the twistor spaces of quaternionic K{\"a}hler manifolds
 of positive scalar curvature (c.~f.~\cite{lebrun2} or \cite{killing}).
 In this sense our estimate improves inequality (\ref{degree}) for
 K\"ahler--Einstein manifolds with complex contact structures.
 \pfill

 Translating the estimate on the degree of the twistor space into
 an estimate on the quaternionic volume we obtain an upper bound
 for the normalized Riemannian volume of a quaternionic K{\"a}hler
 manifold of positive scalar curvature with the same improvement by
 a factor $2$ compared to the estimate given in \cite{sala2}:

 \begin{Corollary}
  Let $(M^{4n},\,g)$ be a quaternionic K{\"a}hler manifold of positive
  scalar curvature. The quaternionic volume of $M$ is bounded
  from above by the quaternionic volume of $\H P^n$:
  $$
   v(M) \;\le\; v(\H P^n) \;=\; 4^n \ .
  $$
 \end{Corollary}

\section{The Hilbert Polynomial of the Wolf Spaces}\label{hwolf}
 The {\it Wolf spaces} are the quaternionic K{\"a}hler symmetric 
 spaces of positive scalar curvature $\kappa>0$. By a classical
 result of Wolf (c.f.~\cite{wolf}) the Wolf spaces correspond up to isometry
 exactly to the simple compact Lie algebras. In particular there are
 three families of Wolf spaces in arbitrary dimensions $4n,\,n\geq 2,$
 namely
 $$
  \H P^n = \frac{\Sp(n+1)}{\Sp(n)\times \Sp(1)},\quad
  \Gr_2(\C^{n+2}) = \frac{\U(n+2)}{\U(n)\times \U(2)},\quad
  \Gr_4(\R^{n+4}) = \frac{\SO(n+4)}{\SO(n)\times \SO(4)}
 $$
 and moreover 5 exceptional Wolf spaces
 $$
  \frac{\mathbf{G}_2}{\SO(4)},\quad
  \frac{\mathbf{F}_4}{\Sp(3)\Sp(1)},\quad
  \frac{\mathbf{E}_6}{\SU(6)\Sp(1)},\quad
  \frac{\mathbf{E}_7}{\mathbf{Spin}(12)\Sp(1)},\quad
  \frac{\mathbf{E}_8}{\mathbf{E}_7\Sp(1)}
 $$
 in dimensions $4n$ with $n=2,7,10,16$ and $28$ respectively.
 It is known that up to isometry there are only finitely many quaternionic 
 K\"ahler manifolds of positive scalar curvature in each dimension  
 and it is natural to conjecture that every quaternionic
 K{\"a}hler manifold $M$ of positive scalar curvature has to be a
 Wolf space. In fact this conjecture has been proved in quaternionic
 dimensions $n=2, 3$  and in dimension $n=4$ under the additional
 assumption $\,b_4(M)=1$. In all known proofs the properties of the
 Hilbert polynomial played a crucial role, providing the main motivation for 
 studying it in detail, in particular we are interested in closed formulas 
 for the Hilbert polynomial of the Wolf spaces.

 Following the general construction of the Wolf spaces given in
 \cite{wolf} let $\g$ be a simple Lie algebra, $\t$ a maximal torus
 and $\l_{\W}\in\t^*$ the highest weight of the adjoint representation
 in a suitable ordering of roots. Let $\<\cdot,\cdot>$ be an invariant
 scalar product on $\g$ and $\g^*$ respectively. The crucial point is
 that the scalar product with $\l_{\W}$  takes exactly five different
 values on the set $\Delta\subset\t^*$ of all roots:
 $$
  \Delta_\hbar\;\;:=\;\;\{\,\mu\in\Delta\,\vert\;\<\l_{\W},\mu>
  \,=\,\hbar\<\l_{\W},\l_{\W}>\,\}
  \qquad\qquad\hbar\;\;=\;\;-1,\,-{1\over2},\,0,\,{1\over2},\,1
 $$
 and the Lie algebra of $\g$ is graded correspondingly into:
 $$
  \g\;\;=\;\;\g_1\,\oplus\,\g_{1\over2}\,\oplus\,
          \g_0\,\oplus\,\g_{-{1\over2}}\,\oplus\,\g_{-1}
  \qquad\qquad
  \g_\hbar\;\;:=\;\;\bigoplus_{\mu\in\Delta_\hbar}\g_\mu
 $$
 In particular the decomposition $\g\,=\,\k\oplus\p\,:=\,(\g_1\oplus
 \g_0\oplus\g_{-1})\,\oplus\,(\g_{1\over2}\oplus\g_{-{1\over2}})$ defines
 a symmetric pair $(\g,\k)$ and the symmetric space corresponding to the
 compact real form of this symmetric pair is the Wolf space associated to
 the simple Lie algebra $\g$. Note that the Lie subalgebra $\k$ decomposes
 further into $\k=\sp(1)_{\W}\oplus\k^0$ with $\sp(1)_{\W}:=\g_1\oplus
 [\g_1,\g_{-1}]\oplus\g_{-1}$ and $\k^0:=\g_0\ominus[\g_1,\g_{-1}]$:

 \begin{Lemma}\label{uniq}
 \summary{Uniqueness of the Minimal Representation}
  For a suitable ordering of roots the Wolf root $\l_{\W}$ is the highest
  weight of the adjoint representation and the set of positive roots contains
  $\Delta^+\supset\Delta_1\cup\Delta_{1\over2}$. Consequently the half sum
  of positive roots is given by
  $$
   \rho\;\;=\;\;{n+1\over2}\l_{\W}\;+\;\rho_{\k^0}
  $$
  where $n$ is the quaternionic dimension of the Wolf space $G/K$ or
  equivalently $n\,=\,\frac14\mathrm{dim}\,\p$.  Moreover for any 
  $r\geq 0$ the finite dimensional representation $\pi=\pi_{r\l_{\W}}$ 
  with highest weight $r\l_{\W}$ is the unique representation $\pi$
  with $\mathrm{Hom}_K(\pi,\,\S^{2r} H)\,\neq\,\{0\}$ and
  $$
   \mathrm{Cas}_{\pi}\;\;=\;\;\l_{2r}\;=\;{\kappa\over2n(n+2)}r(n+1+r)\\
  $$
  In fact $\mathrm{dim}\,\mathrm{Hom}_K(\pi_{r\l_{\W}},\,\S^{2r}H)
  \,=\,1$ and $\pi_{r\l_{\W}}$ occurs once in the sections of $\S^{2r}H$.
 \end{Lemma}

 \proof
 We start by choosing an irrational vector $v_{\k^0}\in\{\l_{\W}\}^\perp$
 orthogonal to $\l_{\W}$ to order the roots of the subalgebra $\k^0$.
 To get an ordering for the roots of $\g$ we make an ansatz $v=\l_{\W}
 +\epsilon v_{\k^0}$ and let $\epsilon$ tend to $0$. Then the scalar product
 of the roots in $\Delta_\hbar$ with $v$ tend to $\hbar\<\l_{\W},
 \l_{\W}>$. For sufficiently small $\epsilon$ the unique root
 $\l_{\W}\in\Delta_1$ will thus be maximal and $\Delta^+\supset
 \Delta_1\cup\Delta_{1\over2}$. Now every root in $\Delta_{1\over2}$
 is of the form $\mu={1\over2}\l_{\W}\,+\,\pr_0(\mu)$ with $\pr_0(\mu)\in
 \{\l_{\W}\}^\perp$ and the half sum of all positive roots is given by:
 \begin{eqnarray*}
  \rho
  &=&{1\over2}\sum_{\mu\in\Delta_1}\mu
     \,+\,{1\over2}\sum_{\mu\in\Delta_{1\over 2}} \mu
     \,+\,{1\over2}\sum_{\mu\in\Delta^+_0}\mu\\
  &=&{1\over2}\l_{\W}\,+\,{2n\over4}\l_{\W}
     \,+\,\sum_{\mu\in\Delta_{1\over2}}\pr_0(\mu)
     \,+\,\rho_{\k^0}
 \;\;=\;\;{n+1\over2}\l_{\W}\,+\,\rho_{\k^0}
 \end{eqnarray*}
 In fact the sum $\sum_{\mu\in\Delta_{1\over2}}\pr_0(\mu)$ must
 vanish, because it is invariant under the Weyl group $W_{\k^0}$
 of $\k^0$. Before we proceed we note the following relation
 between the scalar curvature $\kappa$ of the Wolf space $G/K$
 and the length of the Wolf root:
 \begin{equation}
  {\kappa\over 8n(n+2)}\;\;=\;\;{1\over4}\<\l_{\W},\l_{\W}>
 \end{equation}
 This equation is invariant under rescaling of the scalar product
 $\<\cdot,\cdot>$ on $\g$ and it is thus sufficient to check it in the
 Killing normalization, where the scalar curvature is $\kappa=2n$ and the
 Casimir eigenvalue of the adjoint representation with highest weight
 $\l_{\W}$ is:
 $$
  1\;\;=\;\;\mathrm{Cas}_{\mathrm{ad}}\;\;=\;\;\<\l_{\W},\l_{\W}
  +2\rho>\;\;=\;\;(n+2)\,\<\l_{\W},\l_{\W}>
 $$
 With all these properties of the root system of $\g$ established we
 are now going to show that the representation $\pi_{r\l_{\W}}$ is the
 unique irreducible representation with Casimir eigenvalue $\mathrm{Cas}_{\pi}
 \,=\,r(n+1+r)\<\l_{\W},\l_{\W}>$ occurring in the sections of $\S^{2r}H$.
 Consider an irreducible representation $\pi_\l$ of $\g$ with highest weight
 $\l$ and $\mathrm{Hom}_K(\pi_\l,\,\S^{2r}H)\,\neq\,\{0\}$. This assumption
 implies that there is an element $w\in W_{\g}$ of the Weyl group of $\g$
 such that $\<w\l,\l_{\W}>\geq r\<\l_{\W},\l_{\W}>$. However the Weyl group
 $W_{\k^0}$ of $\k^0$ fixes $\l_{\W}$ and modifying $w$ by elements of
 $W_{\k^0}$ we may assume that $w\l$ is of the form
 $$
  w\l\;\;=\;\;{\delta\over2}\l_{\W}\,+\,\l_{\k^0}
 $$
 with dominant $\k^0$--weight $\l_{\k^0}\geq 0$. In this case we
 have $\delta\geq 2r$ and hence:
 \begin{eqnarray*}
  \<\l,\l+2\rho>
  &\geq&
  \<w\l,w\l+2\rho>\\
  &\geq&
  {\delta(2n+2+\delta)\over4}\<\l_{\W},\l_{\W}>
  \,+\,\<\l_{\k^0},\l_{\k^0}+2\rho_{\k^0}>\\
  &\geq&
  r(n+1+r)\,\<\l_{\W},\l_{\W}>
 \end{eqnarray*}
 If in addition the Casimir eigenvalue of $\pi_{\l}$ is $\l_{2r}\,=\,
 r(n+1+r)\<\l_{\W},\l_{\W}>$ we must have equality everywhere in this chain
 of inequalities, in particular $\delta\,=\,2r$ and $\l_{\k^0}\,=\,0$ as
 $\l_{\k^0}$ is dominant with $\<\l_{\k^0},\l_{\k^0}+2\rho_{\k^0}>\,=\,0$.
 Working out the details of this argument it is easy to verify that
 $\pi_{r\l_{\W}}$ occurs in the sections of $\S^{2r}H$ with multiplicity
 exactly $1$.
 \qed
 \pfill

 The Hilbert polynomial $P(r)$  is given as the index of
 the twisted Dirac operator $D_{\S^{n+2r}H}$ and by (\ref{index}) it is
 the dimension of the eigenspace of $\Delta_{\S^{2r}H}$ for the minimal
 eigenvalue $\l_{2r}$. On symmetric spaces $G/K$ we know that the operator
 $\Delta_\pi$ coincides with the corresponding Casimir operator. Hence
 Frobenius reciprocity leads to:
 $$
  P(r) \;=\; \dim \ker ( \Delta_{\S^{2r}H} \,-\, \l_{2r} ) \;=\;
  \bigoplus_{\pi \in \hat G \atop \Cas_\pi = \l_{2r}}
  \Hom_K(\pi,\,\S^{2r}H) \otimes \dim \pi \ .
 $$
 Combining this decomposition with Lemma \ref{uniq} we obtain an explicit
 formula for the polynomial $P(r)$ on compact Wolf spaces, i.~e.~on all
 symmetric quaternionic K{\"a}hler manifolds:

 \begin{Corollary}
  The Hilbert polynomial on a compact Wolf space is given by:
  $$
   P(r)
   \;\;=\;\;
   \dim (\pi_{\,r\l_{Ad}})
   \;\;=\;\;
   {n+1+2r\over n+1}\,\prod_{\l\in\Delta_{1\over2}}\,
   \left(1+{2r\over 4\<\rho,\l>}\right)
  $$
 \end{Corollary}

 We will close this section in giving some explicit examples for the
 Hilbert polynomial on Wolf spaces. To begin with consider the 8--dimensional
 Wolf space $\G_2/\SO(4)$. The dimension formula for $\G_2$--representations
 implies:
 $$
  P(r) \;\;=\;\; \frac{1}{120}\,(r+2)(3r+5)(2r+3)(3r+4)(r+1) \ .
 $$
 In particular $P$ is of degree 5 with $P(0)=1$ and $P(1) = \dim\G_2 = 14$
 and zeroes in $r=-1$ and $r=-2$. Moreover the leading coefficient is
 $\,\frac{3}{20}$ leading to a quaternionic volume $9$. For the real
 Grassmannians $\,\Gr_4(\R^{n+4})$ we obtain similarly a polynomial
 $$
 P(r) \;=\; \frac{(n+2r)(n+2r+1)(n+2r+2)}{n^2(n+1)(n+2)}\,
 {n+r \choose n-1} {n+r-1 \choose n-1} \ .
 $$
 of degree $2n+1$ with $P(0)=1$ and $P(1)=\dim \SO(n+4)$. Moreover
 we see that the leading coefficient is $\,\frac{8}{n!(n+2)!}$ so
 that the quaternionic volume of the real Grassmannians is given by
 $$
  v(\Gr_4(\R^{n+4}))=\frac{4}{n+2}\,{2n+1\choose n} \ . 
 $$
 Doing the same calculations for the complex Grassmannians $\,\Gr_2(\C^{n+2})$
 we find again a polynomial of degree $2n+1$ with leading coefficient
 $\,\frac{2}{n!(n+1)!}$:
 $$
 P(r) \;=\; \frac{n+2r+1}{n+1}\, {n+r \choose r}^2  \ ,
 \qquad
 v(\Gr_2(\C^{n+2}))={2n+1\choose n} \ .
 $$
\section{Minimal Eigenspaces and Their Prolongations}\label{finito}
 In Section~\ref{estimate} we have seen that the value of the Hilbert
 polynomial $P(r)$ is the dimension of the eigenspace $\E\,=\,\ker
 (\Delta_{\S^{2r}H} - \l_{2r})$ of the operator $\Delta_{\S^{2r}H}$
 acting on section of $\S^{2r}H$, where the constant $\l_{2r}$ is the
 lower bound for its spectrum. In this sense $\E$ is the minimal
 eigenspace of $\Delta_{\S^{2r}H}$. In this section we show that $\E$ is
 the kernel of a first order twistor operator on $\S^{2r}H$ which turns
 out to be a differential operator of finite type. The general theory of
 differential equations of finite type as explained in Section~\ref{killing}
 provides an upper bound on the dimension of $\E$ by the dimension of the
 associated total prolongation. A more direct approach to estimate the
 dimension of $\E$ taking in this section is to define a filtration on
 $\E$ such that the successive filtration quotients embed into the higher
 prolongations. Either way depends on knowing the prolongations explicitly
 and so we will determine them at the end of this section thus proving the
 upper bound on the dimension of $\E$ given in Proposition~\ref{estimateprop}.

 \medskip
 Recall that every representation $\pi$ of $\Sp(1)\cdot\Sp(n)$ gives rise
 to a globally defined vector bundle $\pi(M)$ on every quaternionic K\"ahler
 manifold $M$. By chance however the representations of $\Sp(1)\cdot\Sp(n)$
 are exactly the real representations of $\Sp(1)\times\Sp(n)$ and are thus
 already defined over $\R$. Calculus on $M$ can thus be formulated either
 in terms of real or in terms of complex vector bundles only. The latter
 choice is more convenient as the representations $H$ and $E$ of $\Sp(1)
 \times\Sp(n)$ are both quaternionic and so we will work with complex
 vector bundles exclusively in particular with the complexified tangent
 bundle $TM\otimes_\R\C\,\cong\,H\otimes E$.

 Consequently we will think of the covariant derivative $\nabla\psi$ of a
 section $\psi$ of the complex vector bundle $\S^{2r}H$ as a section of the
 tensor product $(T^*M\otimes_\R\C)\otimes\S^{2r}H$ which is isomorphic to
 $(H\otimes E)\otimes\S^{2r}H\,\cong\,(\S^{2r+1}H\otimes E) \,\oplus\,
 (\S^{2r-1}H\otimes E)$. Projecting $\nabla\psi$ to both of these summands
 in turn defines two natural twistor operators:
 $$
  D^+_d:\;\;
  \Gamma( \S^{2r}H ) \longrightarrow \Gamma(\S^{2r-1}H \otimes E)
  \qquad
  D^+_u: \;\;
  \Gamma( \S^{2r}H ) \longrightarrow \Gamma(\S^{2r+1}H \otimes E) \ .
 $$
 The Weitzenb{\"o}ck formulas of \cite{wgu} characterize the eigenspace
 $\E$ of $\Delta_{\S^{2r}H}$ corresponding to the minimal eigenvalue
 $\lambda_{2r}\,:=\,\frac\kappa{2n(n+2)}\,{r (n+1+r)}$ as the kernel
 of the operator $D^+_u$:

 \begin{Proposition}
  Let $\psi$ be a section of $\S^{2r}H$. Then $\,\Delta_{\S^{2r}H}\,\psi
  \,=\,\lambda_{2r}\psi\,$ if and only if:
  $$
   D^+_u \psi \;=\; 0\,.
  $$
  In this case  $\psi$ satisfies additionally the equation:
  \begin{equation}\label{kill}
   (D^+_d)^*D^+_d\, \psi 
   \;=\;
   \frac{n}{n+1+r}\,\Delta_{\S^{2r}H}\,\psi
   \;=\;
   \frac{\kappa \, r}{2(n+2)}\,\psi \ .
  \end{equation}
 \end{Proposition}

 Note that this proposition is a generalization of a lemma given in
 \cite{sala1}, where it is shown in the case $r\,=\,1$ that the kernel
 of $D^+_u$ is isomorphic to the space $H^0(Z, {\mathcal O}(L))$ and
 that it can be identified with the space of infinitesimal isometries
 on $M$. Indeed in this case it is easy to see that the operator 
 $D^+_d:\;\S^2H\,\longrightarrow\,H\otimes E$ is a multiple of
 the codifferential $d^*:\;\Lambda^2T^*M\,\longrightarrow\,T^*M$
 restricted to $\S^2H\subset\L^2(T^*M\otimes_\R\C)$. Moreover for
 a section $\psi$ in the kernel of $D^+_u$ we have:
 $$
  \Delta(D^+_d\psi) 
  \;=\;
  D^+_d(\Delta_{\S^2H}\,\psi)
  \;=\;
  \frac{\kappa}{2n}\,D^+_d \psi
  \;=\;
  2\, \Ric\,(D^+_d \psi) \ .
 $$
 Hence for a minimal eigensection $\psi$ the coclosed 1--form 
 $D^+_d\psi$ is dual to a Killing vector field and vice versa.
 
 \pfill
 Choose a point $x \in M $ and define the subspace $\E^l_x$ of
 the minimal eigenspace $\E$ by
 $$
  \E^l_x
  \;:=\;
  \{\, \psi \in \E \; | \;
  0\,=\,\psi(x)\,=\,(\nabla \psi) (x)\,=\,\ldots\,=\,(\nabla^l \psi)(x)
  \,=\,0 \;\}\ .
 $$
 Evidently the sequence $\E^l_x,\,l\geq 0,$ of subspaces defines a
 descending filtration of $\E$
 $$
  \E\,:=\,\E^{-1}_x \;\supset\; \E^0_x \;\supset \;\E^1_x\;\supset
  \;\ldots\;\supset\;\E^\infty_x\,:=\,\bigcap_{l\geq 0}\,\E^l_x
 $$
 depending on the point $x\in M$. The intersection $\E^\infty_x$ is the 
 space of minimal eigensections of $\Delta_{\S^{2r}H}$ vanishing in $x$ to
 infinite order. However the operator $\Delta_{\S^{2r}H}$ satisfies
 the strong unique continuation property (c.~f.~\cite{kazdan}) and
 so there is no point on $M$ in which a non--trivial solution $\psi \in \E$
 can possibly vanish to infinite order. Now $\E$ is the eigenspace of
 an elliptic differential operator on a compact manifold $M$ and is
 thus finite dimensional. Consequently
 $$
  \dim \E \;=\; \dim(\,\E^{-1}_x/\E^0_x)\;+\;
  \sum_{l\geq 0}\,\dim(\E^l_x/\E^{l+1}_x) \ .
 $$
 and we will use this equality to prove the upper bound on the dimension of
 $\E$ by estimating the dimension of the successive filtration quotients
 $\E^l_x/\E^{l+1}_x$. In the course of these calculations we will prove
 that the differential equation $D^+_u\psi\,=\,0$ is of finite type,
 i.~e.~that all its higher prolongations vanish everywhere on $M$.
 Indeed we will show that the $l$--th prolongation $\A^{(l)}$ is given as
 $\A^{(l)}\,=\,\S^{2r-l-1}  H\otimes\S^{l+1}E$ and thus vanishes for 
 $l\geq 2r$.  For partial differential equations of finite type
 the assertion $\E^\infty_x\,=\,\{0\}$ for all $x\in M$ is a corollary of
 the construction of a connection on a suitable vector bundle making all
 solutions of the original partial differential equation parallel. For the
 convenience of the reader the general construction of such a connection
 is sketched in Section~\ref{killing}.

 \medskip
 By definition the filtration quotient $\E^{-1}_x/\E^0_x$ embeds
 into $\S^{2r}H$ via $\,[\psi]\longmapsto\psi(x)$ and hence
 $\,\dim\,\E^{-1}_x/\E^0_x\,\leq\,\dim\,\S^{2r}H$. The dimension of
 the higher filtration quotients $\E^l_x/\E^{l+1}_x$ can be estimated
 similarly by embedding them into the higher prolongations $\A^{(l)}_x,
 \,l\geq 0,$ of the symbol $\A^{(0)}$ of the partial differential
 equation $D^+_u\psi\,=\,0$ defined below. The principal symbol
 $\sigma_{D^+_u}:\;(H\otimes E)\otimes\S^{2r}H\longrightarrow
 \S^{2r+1}H\otimes E$ of the twistor operator $D^+_u$ is
 $\Sp(1)\cdot\Sp(n)$--equivariant and thus fixed up to a 
 non--vanishing constant. Redefining $D^+_u$ if necessary we 
 may assume that $\sigma_{D^+_u}$ is given by $\sigma_{D^+_u}
 (h\otimes e \otimes \psi)\,=\,(h\cdot \psi)\otimes e$ for all
 $h\in H,\,e\in E$ and $\psi\in\S^{2r}H$. In terms of the principal
 symbol we can write the partial differential equation:
 $$
  D^+_u\psi\;\;=\;\;\sigma_{D^+_u}(\,\nabla\psi\,)\;\;=\;\;0\,.
 $$
 The general theory introduced in more detail in Section~\ref{killing}
 proceeds by defining
 $$
  \A\;\;:=\;\;(\,\S^{2r}H\,)\;\oplus\;(\,\S^{2r-1}H\otimes E\,)\;\;\subset
  \;\;\S^{\leq 1}(H\otimes E)\otimes\S^{2r}H
 $$
 and the symbol $\A^{(0)}\,:=\,\S^{2r-1}H\otimes E$ of the differential
 equation $D^+_u\psi\,=\,0$ as the kernel of the principal symbol
 $\sigma_{D^+_u}$ of $D^+_u$. This unfortunate clash of nomenclature is
 confusing at first but unavoidable, hopefully the reader unacquainted
 with these concepts will get at least an idea of why different people
 decided to call a subspace and an endomorphism a symbol.
 In order to define the higher prolongations $\A^{(l)},\,l\geq 0,$
 of $\A^{(0)}$ let us introduce the diagonal map
 \begin{equation}\label{comult}
  \Delta:\;\;\S^{k+l}T^*M\;\longrightarrow\;\S^kT^*M\otimes\S^lT^*M
 \end{equation}
 for every $k,\,l\geq 0$ characterized by $\Delta(\frac1{(k+l)!}\xi^{k+l})
 \,=\,\frac1{k!}\xi^k\otimes\frac1{l!}\xi^l$ for every $\xi\in T^*M$.
 The diagonal map is coassociative in the sense that both ways to map
 $\S^{k+l+m}T^*M$ to the tensor product $\S^kT^*M\otimes\S^lT^*M\otimes
 \S^mT^*M$ using the diagonal $\Delta$ result in the same linear map
 sending $\frac1{(k+l+m)!}\xi^{k+l+m}$ to $\frac1{k!}\xi^k\otimes\frac1{l!}
 \xi^l\otimes\frac1{m!}\xi^m$. In the following it will be convenient to
 consider $\S^kT^*M$ as an abstract vector bundle together with a canonical
 embedding into the tensor product bundle $\bigotimes^kT^*M$ obtained by
 iterating the diagonal map as often as possible:
 $$
  \iota:
  \;\;\S^kT^*M\,\;\longrightarrow\;
  {\textstyle\bigotimes^k}T^*M\,,
  \qquad
  \frac1{k!}\,\xi^k\;\longmapsto\;\xi\otimes\ldots\otimes\xi\,.
 $$
 It is straightforward to check the relation $(\iota\otimes\iota)\circ\Delta
 \,=\,\iota$ either directly or using the definition of $\iota$ in terms of
 the diagonal $\Delta$ and coassociativity. With the help of the diagonal map
 $\Delta$ the higher prolongations $\A^{(l)},\,l>0,$ are defined as the
 kernels of the compositions:
 \begin{eqnarray*}
  \S^{l+1}(H\otimes E)\otimes \S^{2r}H
  &\stackrel\Delta\longrightarrow&
  \S^l(H\otimes E)\otimes(H\otimes E)\otimes\S^{2r}H\\
  &\stackrel{\mathrm{id}\otimes\sigma}\longrightarrow&
  \S^l(H\otimes E)\otimes(\,\S^{2r+1}H\otimes E\,)
 \end{eqnarray*}
 In accordance with the general interpretation of the $l$--th prolongation
 we have the lemma:

 \begin{Lemma}\label{symm}
  Let $\psi$ be a minimal eigensection in $\,\E^l_x$, then
  $ \;(\nabla^{l+1}\psi)(x)\,\in\,\A^{(l)}_x$.
 \end{Lemma}

 \proof
 It is a general fact that the iterated covariant derivative
 $\,(\nabla^{l+1} \psi)(x)$ of a section $\psi$ of a vector bundle
 vanishing to order $l$ in a point $x\in M$ is symmetric in all its
 arguments. For $l=1$ this amounts to say that $(\nabla^2_{X,Y}\psi)(x)
 \,-\,(\nabla^2_{Y,X}\psi)(x)\,=\,(R_{X,Y}\psi)(x)\,=\,0$  and
 the general case $l>1$ is verified using essentially the same
 argument with $\nabla^{l-1}\psi$ instead of $\psi$. In our case
 it follows that for a section $\psi\in\E^l_x$ vanishing in $x$
 to order $l$ we we have
 $$
  (\nabla^{l+1}\psi)(x)
  \;\;\in\;\;
  \S^{l+1}(\,H_x\otimes E_x\,)\otimes\S^{2r}H_x
  \;\;\subset\;\;
  {\textstyle\bigotimes^{l+1}}(\,H_x\otimes E_x\,)\otimes\S^{2r}H_x
 $$
 and it remains to prove that $(\nabla^{l+1}\psi)(x)$ is in the
 $l$--th prolongation $\A^{(l)}$ or equivalently is mapped to $0$
 under $(\id\otimes\sigma_{D^+_u})\circ\Delta$. The crucial observation
 is that the diagram
 $$
 \begin{CD}
  \S^{l+1}(\,H_x\otimes E_x\,)
   \otimes\S^{2r}H_x
  & @>{\iota\otimes\id}>> &
  {\textstyle\bigotimes^{l+1}}(\,H_x\otimes E_x\,)
   \otimes\S^{2r}H_x \\
  @V{\Delta}VV && @V{\id}VV \\
  \S^l(\,H_x\otimes E_x\,)\otimes(\,H_x\otimes E_x\,)
   \otimes\S^{2r}H_x
  & @>{\iota\otimes\iota\otimes\id}>> &
  {\textstyle\bigotimes^l}(\,H_x\otimes E_x\,)\otimes(\,H_x\otimes E_x\,)
   \otimes\S^{2r}H_x \\
  @V{\id\otimes\sigma_{D^+_u}}VV && @V{\id\otimes\sigma_{D^+_u}}VV \\
  \S^l(\,H_x\otimes E_x\,)
   \otimes(\,\S^{2r+1}H_x\otimes E_x\,)
  & @>{\iota\otimes\id}>> &
  {\textstyle\bigotimes^l}(\,H_x\otimes E_x\,)
   \otimes(\,\S^{2r+1}H_x\otimes E_x\,)
  \end{CD}
 $$
 commutes as a consequence of the relation $(\iota \otimes \iota)\circ\Delta
 \,=\,\iota$ mentioned above. Now the principal symbol $\sigma_{D^+_u}$ of
 the twistor operator $D^+_u$ is $\Sp(1)\cdot\Sp(n)$--equivariant and so the
 induced bundle homomorphism is parallel. Consequently the iterated covariant
 derivative $(\nabla^{l+1}\psi)(x)\in\bigotimes^{l+1}(\,H_x\otimes E_x\,)
 \otimes\S^{2r}H_x$ of a section $\psi\in\Gamma(\S^{2r}H)$ is mapped to
 $$
  (\id\otimes\sigma_{D^+_u})(\,\nabla^{l+1}\psi\,)(x)
  \;\;=\;\;
  (\id\otimes\sigma_{D^+_u})\;(\,\nabla^l\,\nabla\psi\,)(x)
  \;\;=\;\;
  \nabla^l(\,D^+_u\psi\,)(x)
 $$
 under the composition in the right column. In particular if $\psi\in\E^l_x$
 is a section in the minimal eigenspace and vanishes in $x$ to order $l$ then
 its iterated covariant derivative $(\nabla^{l+1}\psi)(x)$ is not only in the
 image of $\S^{l+1}(\,H_x\otimes E_x\,)\otimes\S^{2r}H_x$, but it is mapped
 to $0$ under the composition in the right column as well. Because all
 horizontal arrows are injective $(\nabla^{l+1}\psi)(x)$ must be in the
 image of the kernel $\A^{(l)}$ of the composition in the left column.
 \qed

 \pfill
 Using this lemma we have an embedding of the successive
 filtration quotients $\E^l_x/\E^{l+1}_x$ into the higher
 prolongations $\A^{(l)}_x$ given by the well--defined map:
 $$
  \E^l_x/\E^{l+1}_x \;\longrightarrow\;\A^{(l)}_x,\qquad
  [\psi]\;\longmapsto\;(\,\nabla^{l+1}\psi\,)(x) \ .
 $$
 Note that this map is injective by construction and hence
 $\dim(\,\E^l_x/\E^{l+1}_x\,)\,\leq\,\dim\A^{(l)}_x$. In the following
 lemma we will determine the higher prolongations and prove the
 isomorphism $\A^{(l)}\,\cong\,\S^{2r-l-1}H\otimes\S^{l+1}E$,
 which eventually completes the proof of Proposition \ref{estimateprop} by:
 $$
  \dim\;\E \;\;\leq\;\; \dim\,\S^{2r}H\;+\;\sum_{l\geq0}\,\dim\;
  \S^{2r-l-1}H\otimes\S^{l+1}E\;\;=\;\;\dim\;\S^{2r}(\,H \oplus E\,)\,.
 $$

 \begin{Lemma}
 The higher prolongations $\A^{(l)},\,l\geq 0,$ of $\A^{(0)}\,=\,
 \S^{2r-1}H\otimes E$ are given by:
 $$
  \A^{(l)}\;\;\cong\;\;\S^{2r-l-1}H\otimes\S^{l+1}E\,.
 $$ 
 \end{Lemma}

 \noindent
 In general knowing all prolongations $\A^{(l)},\,l\geq 0,$ is not quite
 sufficient as we need to know the inclusion maps $\Delta:\;\A^{(l+1)}
 \longrightarrow (H\otimes E)\otimes\A^{(l)}$ as well. In this case
 however these inclusion maps are $\Sp(1)\cdot\Sp(n)$--equivariant
 and thus essentially fixed by the representations $\A^{(l)}$.

 \pfill
 \proof
 The recursive definition of the higher prolongations discussed in more
 detail in Section \ref{killing} lends itself naturally to a proof by
 induction. We will do so by observing that both cases $l=0$  and $l=-1$
 are trivial, if we interpret $\A^{(-1)}$ as the space $\S^{2r}H$.
 According to the general theory described in  Section~\ref{killing}
 there are exact sequences characterizing $\A^{(l+1)}$ for all $l\geq 0$
 as a subspace of $(H\otimes E)\otimes\A^{(l)}$
 $$
  0\;\longrightarrow\;
  \A^{(l+1)}
  \;\stackrel{\Delta}\longrightarrow\;
  (\,H\otimes E\,)\otimes\A^{(l)}
  \;\stackrel{\mathrm{id}\wedge\Delta}\longrightarrow\;
  \L^2(\,H\otimes E\,)\otimes\A^{(l-1)} \ .
 $$
 By induction hypothesis $\A^{(l+1)}$ is
 the kernel of an $\Sp(1)\cdot\Sp(n)$--equivariant map
 $$
  (\,H\otimes E\,)\otimes(\,\S^{2r-l-1}H\otimes \S^{l+1}E\,)
  \;\stackrel{\mathrm{id}\wedge\Delta}\longrightarrow\;
  \L^2(\,H\otimes E\,)\otimes(\,\S^{2r-l}H\otimes\S^lE\,)
 $$
 and we need to know the restriction of the diagonal map $\Delta$
 to the subspace $\S^{2r-l-1}H\otimes\S^{l+1}E\,\subset\,\S^{l+1}(H\otimes E)
 \otimes\S^{2r}H$ to make good use of this description. However the
 restriction is injective and hence fixed up to a non--vanishing constant
 by $\Sp(1)\cdot\Sp(n)$--equivariance alone. As we are only interested in
 the kernel of $\mathrm{id}\wedge\Delta$ we can ignore the constant and
 assume that the restriction of $\Delta$ agrees with
 $$
  \iota:\;\;
  \S^{2r-l-1}H\otimes\S^{l+1}E\;\longrightarrow\;
  (\,H\otimes E\,)\otimes(\,\S^{2r-l}H\otimes\S^lE\,)
 $$
 given by $\iota(\,\psi\otimes\frac1{(l+1)!}e^{l+1}\,)\,:=\,\sum_\alpha\,
 (dh^\flat_\alpha\otimes e)\otimes(h_\alpha\psi\otimes\frac1{l!}e^l)$ for all
 $\psi\in\S^{2r-l-1}H$ and $e\in E$ where $\{h_\alpha\}$ and $\{dh_\alpha\}$
 is a dual pair of bases for $H$ and $H^*$ and $\flat$ is the musical
 isomorphism $H^*\longrightarrow H$ induced by the symplectic form $\sigma$
 on $H$. Moreover the isomorphism
 $$
  \L^2(\,H\otimes E\,)\;\stackrel\cong\longrightarrow\;
  (\S^2H\otimes\L^2E)\oplus\S^2E,\;(\tilde h\otimes\tilde e)\wedge(h\otimes e)
  \;\longmapsto\;(\tilde h\cdot h\otimes\tilde e\wedge e)\oplus
  \sigma(\tilde h, h)\tilde e\cdot e
 $$
 suggests to compose $\,\mathrm{id}\wedge\iota\,$ with the projections onto the
 two summands. The first composition
 $$
  (\,H\otimes E\,)\otimes(\,\S^{2r-l-1}H\otimes\S^{l+1}E\,)
     \;\longrightarrow\;(\,\S^2H\otimes\L^2E\,)\otimes(\,\S^{2r-l}H
     \otimes\S^lE\,)
 $$
 sends $\,(\tilde h\otimes\tilde e)\otimes(\psi\otimes\frac1{(l+1)!}e^{l+1})\,$
 to $\,\sum_\alpha(\tilde h\cdot dh_\alpha^\flat\otimes \tilde e\wedge e)
 \otimes(h_\alpha\psi\otimes\frac1{l!}e^l)$. Hence it is the tensor product
 of the diagonal multiplication $H\otimes\S^{2r-l-1}H\longrightarrow\S^2H
 \otimes\S^{2r-l}H$ with the symplectic form in the $H$--factor with the
 Koszul boundary $E\otimes\S^{l+1}E\longrightarrow\L^2E\otimes\S^lE$ in the
 $E$--factor. The diagonal multiplication with the symplectic form is always
 injective, whereas the kernel of the Koszul boundary is $\,\S^{l+2}E\,$ 
 so that the kernel of the first composition is:
 $$
  (\,H\otimes\S^{2r-l-1}H\,)\otimes(\,\S^{l+2}E\,)
  \;\;\subset\;\;(\,H\otimes E\,)\otimes(\,\S^{2r-l-1}H\otimes\S^{l+1}E\,)\,.
 $$
 Turning to the second composition
 $$
  (\,H\otimes E\,)\otimes(\,\S^{2r-l-1}H\otimes\S^{l+1}E\,)
     \;\longrightarrow\;\S^2E\otimes(\,\S^{2r-l}H\otimes\S^lE\,)
 $$
 which maps $(\tilde h\otimes\tilde e)\otimes(\psi\otimes\frac1{(l+1)!}
 e^{l+1})\,$ to $\,-\,(\tilde e\cdot e)\otimes(\tilde h\cdot\psi\otimes
 \frac1{l!}e^l)$ we observe that it is up to sign the tensor product
 of the multiplication $H\otimes\S^{2r-l-1}H\longrightarrow\S^{2r-l}H$
 with the so called Pl\"ucker differential $E\otimes\S^{l+1}E\longrightarrow
 \S^2E\otimes\S^lE$. The Pl\"ucker differential is injective for $l\geq 1$,
 whereas the kernel of the multiplication is $\S^{2r-l-2}H$ so that for
 $l\geq 1$ at least the kernel of the second composition is the subspace:
 $$
  (\,\S^{2r-l-2}H\,)\otimes(\,E\otimes\S^{l+1}E\,)
  \;\;\subset\;\;(\,H\otimes E\,)\otimes(\,\S^{2r-l-1}H\otimes\S^{l+1}E\,)\,.
 $$
 Consequently for $l\geq 1$ the kernel of $\,\mathrm{id}\wedge\iota\,$ must
 be a subspace of the intersection
 $$
  (\,H\otimes\S^{2r-l-1}H\,)\otimes(\,\S^{l+2}E\,)
  \;\;\cap\;\;
  (\,\S^{2r-l-2}H\,)\otimes(\,E\otimes\S^{l+1}E\,)
 $$
 of kernels of the two compositions of $\,\mathrm{id}\wedge
 \iota\,$ with the projections onto the two summands $\S^2H\otimes\L^2E$ and
 $\S^2E$ of $\L^2(\,H\otimes E\,)$. Evidently this intersection is just the
 subspace $\,\S^{2r-l-2}H\otimes\S^{l+2}E\,\subset\,(\,H\otimes E\,)\otimes(\,
 \S^{2r-l-1}H\otimes\S^{l+1}E\,)$. On the other hand the kernel of
 $\mathrm{id}\wedge\iota$ clearly contains $\S^{2r-l-2}H\otimes\S^{l+2}E$
 and thus agrees with it. The case $l\,=\,0$ requires extra consideration,
 because the kernel of the second composition is strictly larger than
 $\S^{2r-2}H\otimes E\otimes E$. However this difficulty turns out to
 be superficial, because the intersection of the two kernels is still
 given by $\S^{2r-2}H\otimes\S^2E$.
 \qed
\section{Linear Differential Equations of Finite Type}\label{killing}
 Among the partial differential equations the equations of finite type
 form a subclass with a particularly nice description of the set of solutions.
 It turns out that the solutions correspond to the parallel sections for a
 canonical connection on a suitable fibre bundle. Restricting the general
 case somewhat we will only consider a vector bundle $E$ over a manifold $M$
 both endowed with connections and a linear partial differential equation
 $D(\,\psi\,)\,=\,0$ of finite type with a ``parallel'' linear differential
 operator $D$ acting on sections $\psi$ of $E$. In fact the additional 
 assumption of $D$ being ``parallel''
 reduces the technicalities considerably, because the principal symbol of
 $D$ and all its prolongations will be independent of the point of $M$ in
 question.

 In order to make the condition of $D$ being ``parallel'' precise we
 remark that the connections on $E$ and $M$ allow us to define linear
 $k$--th order differential operators $\nabla^k,\,k\geq 0,$ given by
 iterated covariant derivatives as the composition
 $$
  \Gamma(E)\;\stackrel\nabla\longrightarrow\;
  \Gamma(T^*M\otimes E)\;\stackrel\nabla\longrightarrow\;
  \Gamma(\T^2T^*M\otimes E)\;\stackrel\nabla\longrightarrow\;
  \ldots\;\stackrel\nabla\longrightarrow\;\Gamma(\T^kT^*M\otimes E)\;.
 $$
 It is clear that $\nabla^k$ encodes all $k$--th order covariant derivatives
 of a section and thus the operators $\nabla^0,\,\nabla^1,\,\ldots,\,\nabla^k$
 together encode all partial derivatives of a given section in some
 trivialization of $E$. However this information is organized in a rather
 redundant way and for this reason we are well advised to restrict attention
 to the symmetrized iterated covariant derivatives,~i.e. to the differential
 operator
 $$
  \j^k:\;\;\Gamma(E)\;\stackrel{\nabla^k}\longrightarrow\;
  \Gamma(\T^kT^*M\otimes E)\;\stackrel{\frac1{k!}m}\longrightarrow\;
  \Gamma(\S^kT^*M\otimes E)
 $$
 with
 $$
  \j^k_{X_1\cdot\ldots\cdot X_k}\psi
  \;\;:=\;\;
  \frac1{k!}\sum_\tau\;\nabla^k_{X_{\tau(1)},\ldots,X_{\tau(k)}}\psi
 $$
 and $\j^0\psi\;:=\;\psi$ by definition. Its symbol is the identity
 map from $\S^kT^*M\otimes E$ to itself and so the differential operator
 $$
  \j^{\leq k}:\;\;\Gamma(E)\;\longrightarrow\;\Gamma(\S^{\leq k}T^*M\otimes E),
  \quad\psi\;\longmapsto\;\j^0\psi\,\oplus\,\j^1\psi\,\oplus\,\ldots\,
  \oplus\,\j^k\psi
 $$
 of order $k$ from $E$ to the direct sum $\,\S^{\leq k}T^*M\otimes E\,:=\,
 \bigoplus_{l=0}^k\S^lT^*M\otimes E$ has the universal property of a jet
 operator. Namely for every linear differential operator $\mathfrak{D}$ of
 order $k$ from sections of $E$ to sections of a bundle $F$ there is a unique linear map
 $\sigma_{\mathfrak{D}}:\;\S^{\leq k}T^*M\otimes E\,\longrightarrow\,F$
 of vector bundles such that $\mathfrak{D}\,\psi\,=\,\sigma_{\mathfrak{D}}
 (\j^{\leq k}\psi)$ for all sections $\psi$ of $E$. Even for naturally
 defined operators $\mathfrak{D}$ it is somewhat difficult to give explicit
 formulae for these linear maps $\sigma_{\mathfrak{D}}$ and so we need to
 make extensive use of the universal property in order to ensure existence
 of suitable linear maps below.

 Consequently a general linear partial differential equation for sections $\psi$
 of a vector bundle $E$ can be written $P(\,\j^{\leq k}\psi\,)\,=\,0\,$ for
 some linear map $P:\;\S^{\leq k}T^*M\otimes E\longrightarrow F$ of vector
 bundles and it is natural to call an equation of this form parallel if the
 linear map $P$ is. A differential equation which is not parallel is
 presumably unrelated to the affine geometry of $M$ and $E$ and it seems
 better not to use connections and to apply the general language
 of jets instead. Moreover restricting to parallel linear differential equations
 we avoid various technical problems, in particular the principal symbol
 of the differential equation is a parallel subbundle of $\S^kT^*M\otimes E$
 and so all its prolongations are genuine vector bundles over $M$ with
 induced connections.

 Let us defined the total symbol $\A\,:=\,\ker\,P$ as the kernel of the
 linear map $P$. Its principal symbol $\A^{(0)}\,:=\,\ker\,P\cap\S^kT^*M
 \otimes E$ is the kernel of the restriction $\mathrm{res}\,P$ of $P$ to
 the subspace $\S^kT^*M\otimes E$ of polynomials of strict degree $k$. The
 prolongations $\A^{(l)},\,l\geq 0,$ of the symbol $\A^{(0)}$ are then
 defined as the kernels of the composition:
 $$
  \S^{k+l}T^*M\otimes E\;\stackrel\Delta\longrightarrow\;\S^lT^*M\otimes
  \S^kT^*M\otimes E\;\stackrel{\mathrm{id}\otimes\mathrm{res}\,P}
  \longrightarrow\;\S^lT^*M\otimes F\,,
 $$
 where $\Delta$ is the diagonal map as described in (\ref{comult}).
 For a general linear partial differential equation the $\A^{(l)},\,l\geq 0,$
 are families of vector spaces over $M$ but not necessarily vector bundles,
 under the assumption that $P$ is parallel however $\A^{(l)}$ is a parallel
 subbundle of $\S^{k+l}T^*M$ and in particular is equipped with an induced
 connection.

 According to the definition above the prolongation $\A^{(l)}$ is a subspace
 of $\S^{k+l}T^*M\otimes E$. On the other hand the diagonal maps $\Delta:\;
 \S^{k+l+1}T^*M\otimes E\longrightarrow T^*M\otimes\S^{k+l}T^*M\otimes E$
 provide inclusions $\A^{(l+1)}\,\longrightarrow\,T^*M\otimes\A^{(l)}$ for
 all $l\geq 0$. These inclusions can be used to give a recursive definition
 of the higher prolongation $\A^{(l+1)},\,l\geq 0,$ by an exact sequence
 $$
  0\;\longrightarrow\;
  \A^{(l+1)}\;\stackrel\Delta\longrightarrow\;
  T^*M\otimes\A^{(l)}\;\stackrel{\mathrm{id}\wedge\Delta}\longrightarrow\;
  \L^2T^*M\otimes\A^{(l-1)}\,,
 $$
 where the space $\A^{(-1)}$ for $l\,=\,0$ has to be interpreted as
 $\S^{k-1}T^*M\otimes E$. In fact this exact sequence can be extended to
 the right in the obvious way to a complex, the so called Spencer complex
 of the symbol $\A^{(0)}$, and its exactness at $T^*M\otimes\A^{(l)}$ is
 a direct consequence of the coassociativity of the diagonal map $\Delta$
 combined with the fact that the Koszul complex is exact. In order to deal
 with all prolongations at the same time we define the total prolongation
 up to degree $l$ as the vector bundle
 $$
  \A^{\leq l}\;\;:=\;\;\A\,\oplus\,\A^{(1)}\,\oplus\,\ldots\,\oplus\,
  \A^{(l)}\,,
 $$
 note that the first summand is $\A$ and not $\A^{(0)}$. Eventually
 we want to construct a sequence of injective linear maps
 $$
  I^{\leq l}:\;\;\A^{\leq l}\;\longrightarrow\;\S^{\leq k+l}T^*M\otimes E
 $$
 which has the following characteristic property:

 \begin{Lemma}\label{ckt}
  For every solution $\psi\in\Gamma(E)$ of the differential equation
  $P(\,\j^{\leq k}\psi\,)\,=\,0$ and every $l\geq 0$ there exist sections
  $\eta^0\in\Gamma(\A)$ and $\eta^r\in\Gamma(\A^{(r)}),\,l\geq r\geq 1,$
  such that
  \begin{equation}\label{cprop}
   \j^{\leq k+l}\psi\quad=\quad I^{\leq l}(\;\eta^0\oplus\eta^1\oplus
   \ldots\oplus\eta^l\;)
  \end{equation}
  As the maps $I^{\leq l}$ are injective the sections $\eta^0,\,\ldots,
  \,\eta^l$ are uniquely determined by $\psi$.
 \end{Lemma}

 \noindent
 Evidently we can simply choose the first map $I^{\leq 0}$ to be the inclusion
 of $\A\,=\,\A^{\leq 0}$ into $\S^{\leq k}T^*M\otimes E$ and $\eta^0\,:=\,
 \j^{\leq k}\psi$. We postpone the proof of the lemma for general $l>0$
 for a moment and give a recursive definition of the higher order maps
 $I^{\leq l},\,l>0,$ first. For this purpose we consider the composition
 $\j^l\circ P\circ\j^{\leq k}$, which is a differential operator of order
 $k+l$ from $E$ to $\S^lT^*M\otimes F$. The universal property of jet
 operators ensures the existence of a linear map
 $$
  P^l:\;\;\S^{\leq k+l}T^*M\otimes E\;\longrightarrow\;\S^lT^*M\otimes F
 $$
 such that $P^l(\,\j^{\leq k+l}\psi\,)\,=\,\j^l(\,P(\,\j^{\leq k}\psi\,)\,)$
 for all sections $\psi$ of $E$ and symbolic calculus asserts that this map
 $P^l$ restricted to the subspace $\S^{k+l}T^*M\otimes E$ is simply the
 composition
 $$
  \S^{k+l}T^*M\otimes E\;\stackrel\Delta\longrightarrow\;\S^lT^*M\otimes
  \S^kT^*M\otimes E\;\stackrel{\mathrm{id}\otimes\mathrm{res}\,P}
  \longrightarrow\;\S^lT^*M\otimes F
 $$
 with kernel $\A^{(l)}$.  The crucial step in the construction is
 now the choice of a partial inverse
 $$
  S^l:\;\S^lT^*M\otimes F\;\longrightarrow\;\S^{k+l}T^*M\otimes E
 $$
 for the restriction of $P^l$ to $\S^{k+l}T^*M\otimes E$ in the sense that
 $S^l$ maps every element in the image of $\S^{k+l}T^*M\otimes E$ under $P^l$
 to some preimage, in other words $P^l\circ S^l\circ P^l\,=\,P^l$ holds true
 on the subspace $\S^{k+l}T^*M\otimes E$. Clearly such a partial inverse
 exists, because the kernel $\A^{(l)}$ of this restriction is a vector
 bundle on $M$ and so is its image. With the inverses $S^l,\,l\geq 0,$
 chosen we can inductively define the sequence the higher order maps by
 setting
 $$
  I^{\leq l+1}(\;\eta^0\,\oplus\,\ldots\,\oplus\eta^{l+1}\;)
  \;\,:=\;\,I^{\leq l}(\;\eta^0\oplus\ldots\oplus\eta^l\;)\oplus
  \left(\,\eta^{l+1}\,-\,S^{l+1}\,P^{l+1}\,(\,I^{\leq l}(\;
  \eta^0\oplus\ldots\oplus\eta^l\;)\oplus 0\;)\,\right)
 $$
 for $l\geq 0$ with $I^{\leq 0}$ being the inclusion as above. Obviously
 $\A^{(r)},\,r\geq 1$ is mapped to the direct sum of the $\S^{k+l}T^*M
 \otimes E$ with $l\geq r$ only and the induced map $\A^{(r)}\longrightarrow
 \S^{k+r}T^*M\otimes E$ is simply the inclusion. Hence it follows by
 straightforward induction that all $I^{\leq l}$ are injective. Moreover
 the maps $I^{\leq l},l\ge 0$, have the characteristic property  claimed
 in Lemma~\ref{ckt}:
 
 \medskip
 \proof
 Starting the induction with $l=0$ we simply
 choose $\eta^0\,:=\,\j^{\leq k}\psi$ to be the $k$--jet of $\psi$ in
 $\A\,=\,\A^{\leq 0}$ as already discussed above. Consider now a solution
 $\psi\in\Gamma(E)$ and a fixed point of $M$. By induction hypothesis we
 may assume that in this point we have equality
 $$
  \j^{\leq k+l+1}\psi
  \quad=\quad I^{\leq l}(\;\eta^0\oplus\ldots\oplus\eta^l\;)\;\oplus\;
  \j^{k+l+1}\psi
 $$
 for suitably chosen $\eta^0\oplus\eta^1\oplus\ldots\oplus\eta^l$ in
 $\A^{\leq l}$. According to our assumption $\psi$ is a solution to
 $P(\,\j^{\leq k}\psi\,)\,=\,0$ and using the fact the $P^{l+1}$ is
 linear we get:
 \begin{eqnarray*}
  P^{l+1}(\;\j^{\leq k+l+1}\psi\,)
  &=&
  \j^{l+1}\,(\,P(\,\j^{\leq k}\psi\,)\,)\;\;=\;\; 0\\
  &=&
  P^{l+1}(\;I^{\leq l}(\eta^0\oplus\ldots\oplus\eta^l)\oplus 0\;)
  \;+\;P^{l+1}(\;0\oplus\j^{k+l+1}\psi\;)\,.
 \end{eqnarray*}
 In particular $P^{l+1}(\,I^{\leq l}(\eta^0\oplus\ldots\oplus\eta^l)
 \oplus 0\,)$ lies in the image of $\S^{k+l+1}T^*M\otimes E$ under $P^{l+1}$
 and as $S^{l+1}$ is a partial inverse for $P^{l+1}$ on this image we may
 write this equality as:
 $$
  P^{l+1}(\;\j^{k+l+1}\psi\;)\;+\;P^{l+1}\,S^{l+1}\,P^{l+1}
  (\;I^{\leq l}(\eta^0\oplus\ldots\oplus\eta^l)\oplus 0\;)\quad=\quad0\,.
 $$
 Hence the two elements $\j^{k+l+1}\psi$ and $S^{l+1}P^{l+1}(\,I^{\leq l}
 (\eta^0\oplus\ldots\oplus\eta^l)\oplus 0\,)$ of $\S^{k+l+1}T^*M\otimes E$
 only differ by an element $\eta^{l+1}$ in the kernel $\A^{(l+1)}$ of
 $P^{l+1}$ restricted to $\S^{k+l+1}T^*M\otimes E$. In other words
 $\j^{k+l+1}\psi\;=\;\eta^{l+1}\,-\,S^{l+1}P^{l+1}(\,I^{\leq l}(\eta^0
 \oplus\ldots\oplus\eta^l)\oplus0\,)$ and we conclude
 \begin{eqnarray*}
  \j^{\leq k+l+1}\psi
  &=& \j^{\leq k+l}\psi\oplus\j^{k+l+1}\psi\\
  &=& I^{\leq l}(\;\eta^0\oplus\ldots\oplus\eta^l\;)
   \oplus\left(\,\eta^{l+1}\,-\,S^{l+1}P^{l+1}(\,I^{\leq l}
   (\;\eta^0\oplus\ldots\oplus\eta^l\;)\oplus 0\,)\,\right)\\
  &=& I^{\leq l+1}(\;\eta^0\oplus\ldots\oplus\eta^l\oplus\eta^{l+1}\;)
 \end{eqnarray*}
 using again the induction hypothesis. It is clear that the section
 $\eta^0\oplus\ldots\oplus\eta^l\oplus\eta^{l+1}$ of $\A^{\leq l+1}$
 constructed pointwise this way is a smooth section of $\A^{\leq l+1}$.
 \qed

 \begin{Definition}
 A partial differential equation $P(\,\j^{\leq k}\psi\,)\,=\,0$ is called
 of finite type if and only if its higher prolongations $\A^{(l)}$
 vanish everywhere on $M$ for all sufficiently large $l\gg 0$. According to
 the recursive definition of the prolongations $\A^{(l)},\,l>0,$ a partial
 differential equation is of finite type if and only if $\A^{(d+1)}\,=\,
 \{0\}\,$ everywhere for some $d\geq 0$.
 \end{Definition}

 Differential equations of finite type have a very neat characterization
 in terms of their symbol. Recall that a differential operator is called
 elliptic if and only if its principal symbol considered as a linear map
 is injective for every real cotangent vector. Thinking of the principal
 symbol as a subspace this means that the principal symbol $\A^{(0)}$ of
 an elliptic differential operator of order $k$ does not contain any
 $\frac1{k!}\,\xi^k\otimes e\neq 0$ with $\xi\in T^*M$ and $e\in E$:

 \begin{Theorem}
 A differential equation is of finite type if and only
 if the complexification
 $$
  \A^{(0)}\otimes_\R\C\;\;\subset\;\;(\,\S^kT^*M\otimes E\,)\otimes_\R\C
  \;\;\cong\;\;\S^k(\,T^*M\otimes_\R\C\,)\,\otimes\,(\,E\otimes_\R\C\,)
 $$
 of its principal symbol does not contain any $\frac1{k!}\,\xi^k\otimes e
 \neq 0$ with $\xi\in T^*M\otimes_\R\C$ and $e\in E\otimes_\R\C$.
 \end{Theorem}

 \pfill
 Using this criterion one can easily check that the operator $D^+_u$
 of Section~5 is of finite type. Other examples are the classical twistor
 operator in spin geometry or the twistor operator defining conformal
 Killing forms (c.~f.~\cite{confkill}). Explicit upper bounds for the
 dimension of the solution space of a differential equation of finite
 type however can only be found by calculating all non--vanishing
 prolongations.
 
 \pfill
 Of course Lemma \ref{ckt} is valid for general parallel linear
 differential equations $P(\,\j^{\leq k}\psi\,)\,=\,0$. What is remarkable
 about differential equations of finite type however is that the sequence
 $\A^{(l)},\,l\geq 0,$ of prolongations of the symbol becomes trivial
 $\A^{(l)}\,=\,\{0\}$ for sufficiently large $l > d\geq 0$. From
 this point on we have isomorphisms
 $$
  \A^{\leq d}\;\;=\;\;\A^{\leq d+1}\;\;=\;\;\A^{\leq d+2}\;\;=\;\;
  \ldots\;\;=\;\;\A^{\leq l}\,,
 $$
 but we may still construct the maps $I^{\leq l}$ for all $l\geq 0$ as above.
 Thus for every $l\geq d$ the $k+l$--jet of some solution $\psi$ at a point of
 $M$ is determined by the values of the sections $\eta^0,\,\ldots,\,\eta^d$,
 which in turn are determined by the $k+d$--jet of $\psi$. In essence this
 means that from the point $l\,=\,d$ on we can calculate all higher order jets
 of a solution $\psi$ in a point from the $k+d$--jet of $\psi$ in this point
 alone. This property is very similar to the behavior of parallel sections
 of a vector bundle whose higher order jets at a point are all determined by
 the value in this point.

 In order to make this analogy precise we apply the universal property of
 the jet operators to the differential operator $\nabla\,\j^{\leq k+d}$ of
 order $k+d+1$ from $E$ to $T^*M\otimes\S^{\leq k+d}T^*M\otimes E$ and find
 linear maps
 $$
  \widehat{\Delta}:
  \;\;\S^{\leq k+d+1}T^*M\otimes E\;\longrightarrow\;T^*M\otimes
  \S^{\leq k+d}T^*M\otimes E
 $$
 such that $\nabla(\,\j^{\leq k+d}\psi\,)\,=\,\widehat{\Delta}
 (\,\j^{\leq k+d+1}\psi\,)$ for every section $\psi$ of $E$. Now every
 solution $\psi$ of the differential equation $P(\,\j^{\leq k}\psi\,)\,=\,0$
 corresponds to a section $\eta^0\oplus\ldots\oplus\eta^d$ of the bundle
 $\A^{\leq d}$. Recalling that $\A^{\leq d}$ is equipped with a natural
 connection we conclude
 \begin{eqnarray*}
  \nabla(\,I^{\leq d}(\;\eta^0\,\oplus\,\ldots\,\oplus\,\eta^d\;)\,)
  &=&
  (\,\nabla\,I^{\leq d}\,)(\;\eta^0\,\oplus\,\ldots\,\oplus\,\eta^d\;)
  \,+\,(\mathrm{id}\otimes I^{\leq d})
  (\,\nabla(\eta^0\,\oplus\,\ldots\,\oplus\,\eta^d)\,)\\
  &=&
  \nabla(\,\j^{\leq k+d}\psi\,)\;\;=\;\;
  \widehat{\Delta}
  (\,I^{\leq d+1}(\;\eta^0\,\oplus\,\ldots\,\oplus\,\eta^d\;)\,)
 \end{eqnarray*}
 and so the covariant derivative of the section $\eta^0\oplus\ldots
 \oplus\eta^d$ is completely determined by
 $$
  (\mathrm{id}\otimes I^{\leq d})
  (\;\nabla(\eta^0\,\oplus\,\ldots\,\oplus\,\eta^d)\;)
  \;\;=\;\;
  \widehat{\Delta}
  (\,I^{\leq d+1}(\;\eta^0\,\oplus\,\ldots\,\oplus\,\eta^d\;)\,)
  \;-\;(\,\nabla\,I^{\leq d}\,)(\;\eta^0\,\oplus\,\ldots\,\oplus\,\eta^d\;)
 $$
 as $\,\mathrm{id}\otimes I^{\leq d}$ is injective. We may reformulate
 this property in terms of the differential operator
 \begin{eqnarray*}
  D(\;\eta^0\oplus\ldots\oplus\eta^d\;)
  &:=&
  (\mathrm{id}\otimes I^{\leq d})
  (\;\nabla(\eta^0\,\oplus\,\ldots\,\oplus\,\eta^d)\;)\\
  &&\quad-\;\left(
  \,\widehat{\Delta}
  (\,I^{\leq d+1}(\;\eta^0\,\oplus\,\ldots\,\oplus\,\eta^d\;)\,)
  \;-\;(\,\nabla\,I^{\leq d}\,)(\;\eta^0\,\oplus\,\ldots\,\oplus\,\eta^d\;)\,
  \right)
 \end{eqnarray*}
 from sections of $\A^{\leq d}$ to $T^*M\otimes\S^{\leq k+d}T^*M\otimes E$,
 whose principal symbol is the inclusion:
 $$
  (\mathrm{id}\otimes I^{\leq d}):\;\;
  T^*M\otimes\A^{\leq d}\;\longrightarrow\;
  T^*M\otimes\S^{\leq k+d}T^*M\otimes E\,.
 $$
 Every section $\eta^0\oplus\ldots\oplus\eta^d$ of $\A^{\leq d}$ which
 corresponds to a solution $\psi$ of the differential equation $P(\,
 \j^{\leq k}\psi\,)\,=\,0$ is killed by $D$. More precisely the solutions
 $\psi$ of the original differential equation correspond bijectively to
 the sections of $\A^{\leq d}$ in the kernel of the operator $D$. In fact
 the differential operator $D$ is essentially the restriction of the so
 called canonical connection
 $$
  \nabla^{\mathrm{can}}:\;\;\Gamma(\,\S^{\leq k+d+1}T^*M\otimes E\,)
  \;\longrightarrow\;\Gamma(\,T^*M\otimes\S^{\leq k+d}T^*M\otimes E\,)
 $$
 or the ``naive'' Spencer operator to the image of $\A^{\leq d}$ under
 $I^{\leq d+1}$, and it is well--known and easy to prove that the sections
 of $\S^{\leq k+d+1}T^*M\otimes E$ parallel under the canonical connection
 $\nabla^{\mathrm{can}}$ are precisely the $k+d+1$--jets $\j^{k+d+1}\psi$
 of sections $\psi$ of $E$.


 \pfill
 Now we have succeeded in  identifying the space of solutions to the original
 differential equation $P(\,\j^{\leq k}\psi\,)\,=\,0$ with the kernel of the
 differential operator $D$ and one may wonder what we have achieved at all.
 However the differential operator $D$ has injective symbol and thus the
 condition $D(\,\eta^0\oplus\ldots\oplus\eta^d\,)\,=\,0$ is stronger than
 $\eta^0\oplus\ldots\oplus\eta^d$ being parallel for a suitable connection
 on $\A^{\leq d}$. Indeed for a first order differential operator
 ${\mathfrak D}:\;\Gamma(\A)\,\longrightarrow\,\Gamma(F)$ with injective
 symbol $\sigma_{\mathfrak{D}}$ the image of $\sigma_{\mathfrak{D}}$ is a
 subbundle of $F$ and hence we can choose a smooth projection $\pr$ onto
 this image. The composition $\sigma^{-1}_{\mathfrak{D}}\circ\pr\circ
 \mathfrak{D}$ is a differential operator of first order from $\A$ to
 $T^*M\otimes\A$ with principal symbol given by the identity of
 $T^*M\otimes\A$. In other words $\sigma^{-1}_{\mathfrak{D}}\circ
 \pr\circ\mathfrak{D}$ is a connection on $\A$ making sections in
 the kernel of $\mathfrak{D}$ parallel. In our case we have to compose
 $D$ with some linear projection
 $$
  \pr:\;\;T^*M\otimes\S^{\leq k+d}T^*M\otimes E\;\longrightarrow
  \;(\,\mathrm{id}\otimes I^{\leq d}\,)(\,T^*M\otimes\A^{\leq d}\,)
 $$
 onto the image of $T^*M\otimes\A^{\leq d}$ under $\mathrm{id}\otimes
 I^{\leq d}$ in order to obtain a connection on the vector bundle
 $\A^{\leq d}$ such that every solution to the original equation
 corresponds to a parallel section of $\A^{\leq d}$. However the
 nice feature of the operator $D$ of classifying the solutions
 completely may get lost in projecting $D$ to this Killing connection,
 in other words there may be sections of $\A^{\leq d}$ parallel for
 the Killing connection but not killed by $D$, which do not correspond
 to any solution of the original equation at all:

 \begin{Proposition}
  Let $P(\,\j^{\leq k}\psi\,)\,=\,0$ be a parallel linear differential
  equation of finite type with $\A^{(d+1)}\,=\,\{0\}$ on a connected
  manifold $M$. The dimension of the space of solutions of this differential
  equation is bounded above by the dimension of $\A^{\leq d}$. Moreover the
  solutions of the differential equation correspond bijectively to sections
  of the vector bundle $\A^{\leq d}$ parallel for every connection on
  $\A^{\leq d}$ constructed from projecting the canonical connection
  $D$ as above.
 \end{Proposition}

\begin{thebibliography}{MMM99-9}
%
\bibitem[Ka88]{kazdan}
 \textit{Kazdan, J.}
 Unique continuation in geometry. 
 Comm. Pure Appl. Math. {\bf 41} (1988), no. 5, 667--681. 

\bibitem[KSW98]{wgu}
 \textit{Kramer, W., Semmelmann, U. and Weingart, G.}
 Quaternionic Killing spinors.
 Ann. Glob. Anal. Geom. {\bf 16} (1998), 63-87. 

\bibitem[LS94]{lebrun1}
  \textit{LeBrun, C. and Salamon, S.}
  Strong rigidity of positive quaternion--K{\"a}hler manifolds.
  Invent. Math. {\bf 118} (1994), no.1, 109--132.

\bibitem[Le95]{lebrun2}
  \textit{LeBrun, C.}
  Fano manifolds, contact structures, and quaternionic geometry.
  Internat. J. Math. {\bf 6} (1995), no. 3, 419--437.

\bibitem[MS96]{killing}
   \textit{Moroianu, A. and Semmelmann, U.}
   K{\"a}hlerian {K}illing spinors, complex contact structures and
    twistor spaces. 
   C. R. Acad. Sci. Paris Ser. I Math. {\bf 323} (1996), no. 1, 57--61. 

\bibitem[Sa82]{sala1}
  \textit{Salamon, S.M.}
  Quaternionic {K}{\"a}hler manifolds,
  Invent. Math. {\bf 67} (1982), 143--171.

\bibitem[Sa99]{sala2}
  \textit{Salamon, S.M.}
  Quaternion--{K}{\"a}hler geometry.
  Surveys in differential geometry: essays on Einstein manifolds, 
  83--121, Surv. Differ. Geom., VI, Int. Press, Boston, MA, (1999). 

\bibitem[SW02]{betti}
  \textit{Semmelmann, U. and Weingart, G.},
  Vanishing theorems for quaternionic {K}{\"a}hler manifolds ,
  J. Reine Angew. Math. {\bf 544} (2002), 111--132.

\bibitem[Se02]{confkill}
  \textit{Semmelmann, U.},
  Conformal Killing forms on Riemannian manifolds,
  preprint (2002), math.DG/0206117.

\bibitem[Wo65]{wolf}
  \textit{Wolf, J.},
  Complex homogeneous contact manifolds and quaternionic
   symmetric spaces,
  J. Math. Mech. {\bf 14} (1965), 1033--1047;
%
\end{thebibliography}
\end{document}